\documentclass[11pt]{article}

\usepackage{amsmath,amssymb,amsthm}
\usepackage{hyperref}
\usepackage{float}

\DeclareMathOperator{\GL}{GL}
\DeclareMathOperator{\tr}{tr}
\DeclareMathOperator{\Irr}{Irr}
\DeclareMathOperator{\cd}{cd}
\DeclareMathOperator{\dl}{dl}
\DeclareMathOperator{\Sub}{Sub}

\newtheorem{theorem}{Theorem}[section]
\newtheorem{conjecture}[theorem]{Conjecture}
\newtheorem{lemma}[theorem]{Lemma}
\newtheorem{proposition}[theorem]{Proposition}

\title{The Main Problem of Block Theory:\\
Picky Elements and Subnormalizers}

\author{Alexander Moret\'o\\
Department of Mathematics\\
Universitat de Val\`encia\\
\texttt{alexander.moreto@uv.es}}

\date{}

\begin{document}

\maketitle

\begin{abstract}
This article is essentially an English translation of my paper
\cite{gaceta}, published in \emph{La Gaceta de la RSME}. Its aim is to
present, for a broad mathematical audience, a research programme in local
representation theory that goes beyond the classical restrictions to
characters of $p'$-degree, characters of height zero, and blocks of abelian
defect. The final and most recent part of this programme concerns
Alperin's main problem of block theory: the search for local rules for
character values. In that direction I describe the conjectures on picky
elements and subnormalizers, which suggest that the sets $\Irr^x(G)$ and
the subgroups $\Sub_G(x)$ are the natural objects attached to a
$p$-element $x$.
\end{abstract}

\section{Introduction}
\label{sec:intro}

Finite group theory experienced a major expansion in the 1960s and 1970s,
driven by the belief that the classification of the finite simple groups
could be achieved. Recall that a finite group $G$ is called simple if its
only normal subgroups are $\{1\}$ and $G$ itself. A fundamental ingredient
in reaching this goal was the work of Richard Brauer on modular
representation theory, which opened entirely new directions in the subject.

After the announcement of the classification, around 1980, it was sometimes
suggested that group theory was essentially finished, in the sense that no
fundamental problems remained to be solved. A reflection of this intellectual
climate was the conference held in Santa Cruz in 1979. In the preface to the
proceedings, Geoffrey Mason wrote:

\medskip

\noindent
\emph{``In the last year or so there have been widespread rumors that group
theory is finished, that there is nothing more to be done. It is not so.
While it is true that we are tantalizingly close to that pinnacle
representing the classification of finite simple groups, one should remember
that only by reaching the top can one properly look back and survey the
neighboring territory.''}

\medskip

Time has shown that Mason was right. In particular, in his contribution to
the Santa Cruz proceedings, Jon Alperin formulated a number of conjectures,
both his own and those of other authors, that have profoundly shaped the
development of representation theory ever since. Many of them have been
resolved only recently; others appear close to resolution.

This might suggest, once again, that the representation theory of finite
groups is nearly complete. However, as I shall explain in this article, the
solution of these classical conjectures does not exhaust the deeper structure
of blocks. On the contrary, recent advances point toward a much broader
picture, going well beyond the traditionally studied cases, such as groups
with abelian Sylow $p$-subgroups, blocks of abelian defect, or characters of
degree not divisible by a prime $p$.

My aim is to present some of these classical problems and, above all, to
explain why more general approaches suggest a much larger neighboring
territory. In particular, I argue that several fundamental questions in
block theory appear to have natural extensions to arbitrary defect,
beyond the abelian setting in which they were traditionally studied. From this perspective, the recent successes do not mark an endpoint,
but rather the beginning of an exploration that has only just begun.

In Section~\ref{sec:prel} I introduce the basic notions on finite groups,
representations, and blocks that will appear throughout the paper. In
Section~\ref{sec:local} I recall some of the main problems in representation
theory from 1980 to the present. In Section~\ref{sec:future} I state
extensions of these problems beyond the abelian or height-zero setting that
I proposed in earlier work. The central part of the article is
Section~\ref{sec:picky}, where I present several conjectures related to
what Alperin called the ``main problem of block theory.'' A version of these
ideas addressed to specialists, including more general formulations and
further results, will be posted on the arXiv in the near future
\cite{mor}. These conjectures are likely to play an important role in
future work.

\section{Preliminaries}
\label{sec:prel}

In this section I briefly review some basic notions concerning finite
groups, representations, and blocks. The goal is not to give a
comprehensive introduction, but simply to fix the language needed later on.

\subsection{Finite groups}
\label{sub:groups}

Throughout this article all groups are finite. The number of elements of a
group $G$ is called its order and is denoted by $|G|$.

Two elements $x,y \in G$ are said to be conjugate if there exists $g \in G$
such that $y=g^{-1}xg$. Conjugacy classes play a fundamental role in the
theory.

We write
\[
G'=[G,G]=\langle [x,y]\mid x,y\in G\rangle,
\]
where $[x,y]=x^{-1}y^{-1}xy$ is the commutator. The group $G$ is abelian if
and only if $G'=1$.

We define the derived series
\[
G \geq G' \geq G'' \geq G''' \geq \cdots,
\]
where $G''=(G')'$, $G'''=(G'')'$, and so on. The group $G$ is called
solvable if this series reaches the trivial subgroup. In that case, the
derived length $\dl(G)$ is the smallest integer $n$ such that the $n$-th
derived subgroup is trivial. In particular, $G$ is abelian if and only if
$\dl(G)=1$.

Let $p$ be a prime. A $p$-subgroup of $G$ is a subgroup whose order is a
power of $p$. If
\[
|G|=p^a m, \qquad p \nmid m,
\]
then a Sylow $p$-subgroup of $G$ is a subgroup of order $p^a$. By Sylow's
theorems, such subgroups always exist and are all conjugate. In particular,
any two Sylow $p$-subgroups are isomorphic. Moreover, every $p$-subgroup of
$G$ is contained in a Sylow $p$-subgroup.

These subgroups play a central role in group theory and representation
theory. For instance, they will appear as defect groups of principal
blocks.

A finite group $G$ is called nilpotent if for every prime $p$ it has a
unique Sylow $p$-subgroup, or equivalently, if it is the direct product of
its Sylow subgroups. Every finite group $G$ has a largest normal nilpotent
subgroup, called the Fitting subgroup and denoted by $F(G)$.

Given an integer $n$ and a prime $p$, I denote by $n_p$ the $p$-part of
$n$, that is, the largest power of $p$ dividing $n$. Throughout the paper,
$G$ will denote a finite group and $p$ a prime.

\subsection{Representations and characters}
\label{sub:chars}

A complex representation of a finite group $G$ is a homomorphism
\[
\rho : G \to \GL(n,\mathbb{C}),
\]
where $\GL(n,\mathbb{C})$ is the group of invertible $n\times n$
complex matrices. The associated character is the function
$\chi : G \to \mathbb{C}$ defined by
\[
\chi(g)=\tr(\rho(g)),
\]
where $\tr$ denotes the trace.

Characters are constant on conjugacy classes. A character is called
irreducible if it cannot be written as a sum of other characters. I denote
by $\Irr(G)$ the set of irreducible characters of $G$.

A basic result states that
\[
|\Irr(G)| = k(G),
\]
where $k(G)$ is the number of conjugacy classes of $G$. Thus, the values of
the irreducible characters on the conjugacy classes form a square matrix of
size $k(G)$, known as the character table of $G$.

The constant function $1_G(g)=1$ for all $g\in G$ is called the principal
character.

For a general introduction to character theory I refer to the classical
book of Isaacs.

\subsection{Brauer blocks}
\label{sub:blocks}

If we work over an algebraically closed field of characteristic $p>0$
instead of over $\mathbb{C}$, we obtain the modular representation theory
of $G$. This theory was developed by Brauer in order to better understand
complex representations and their relation to the structure of the group.

Given a prime $p$, Brauer's work partitions $\Irr(G)$ into the
so-called $p$-blocks. There are several equivalent descriptions of this
partition; one convenient description for my purposes is the following:
the $p$-blocks are the connected components of the graph whose vertices are
the irreducible characters, where two characters $\chi$ and $\psi$ are
joined if
\[
\sum_{g\in G_{p'}} \chi(g)\psi(g^{-1}) \neq 0,
\]
where $G_{p'}$ denotes the set of $p$-regular elements of $G$.

To each block $B$, Brauer associated a conjugacy class of $p$-subgroups of
$G$, called the defect groups of $B$. In a rough sense, the defect group measures the ``complexity'' of the
block: for instance, a block has
trivial defect group if and only if it contains exactly one irreducible
character.

The block containing the principal character is called the principal block.
Its defect groups are precisely the Sylow $p$-subgroups of $G$.

For more details on the classical theory of Brauer blocks I refer to the
literature.

\section{Local-global problems}
\label{sec:local}

As mentioned in the introduction, Brauer's work was fundamental for the
classification of the finite simple groups. Although the classification
was announced in 1983, the final part was not published until 2004, in the
work of Aschbacher and Smith.

The classification increased the hope that the so-called local-global
problems in the representation theory of finite groups could be solved.
Since then, a large part of the research in the area has focused on these
questions. These problems have guided not only the solution of specific
questions, but also the way in which block theory has been understood.

There is no precise formal definition of a ``local subgroup,'' but the idea
is that such subgroups are closely related to the $p$-local structure of
$G$, such as normalizers of $p$-subgroups or centralizers of $p$-elements.
The prototypical example is, of course, a Sylow subgroup.

\subsection{The It\^o--Michler theorem}
\label{sub:ito}

The first major local-global result obtained using the classification of
the finite simple groups was the It\^o--Michler theorem.

\begin{theorem}[It\^o--Michler]
\label{thm:ito}
Let $G$ be a finite group and $p$ a prime. Then $G$ has a normal abelian
Sylow $p$-subgroup if and only if $p$ does not divide $\chi(1)$ for all
$\chi \in \Irr(G)$.
\end{theorem}

The numbers $\chi(1)$, called the degrees of irreducible characters,
reflect the size of the matrices appearing in the associated
representations. The theorem relates a global property (no character
degree is divisible by $p$) with a local property (the Sylow $p$-subgroup
is normal and abelian).

The difficult part of the proof is to show that if $p$ does not divide
$\chi(1)$ for any $\chi \in \Irr(G)$, then the Sylow $p$-subgroup
is normal. It\^o reduced this problem in the 1950s to the case where $G$
is simple, and Michler completed the proof almost thirty years later
using the classification.

\subsection{Brauer's height zero conjecture}
\label{sub:bhz}

Brauer's height zero conjecture can be viewed as a block version of the
It\^o--Michler theorem. Brauer proved that if $B$ is a $p$-block of $G$
with defect group $D$, $|D|=p^d$, then any irreducible character $\chi$
in $B$ satisfies
\[
\chi(1)_p = p^{a-d+h}
\]
for some integer $h \geq 0$, called the height of $\chi$. The conjecture,
formulated in 1955, states the following.

\begin{conjecture}[Brauer's height zero conjecture]
\label{conj:bhz}
Let $G$ be a finite group, $p$ a prime, and $B$ a $p$-block of $G$ with
defect group $D$. Then $D$ is abelian if and only if all irreducible
characters in $B$ have height zero.
\end{conjecture}

Its proof, obtained in \cite{km,mnst,ruh},
is one of the major recent achievements in the area.

In the case of the principal block, the defect group is a Sylow
$p$-subgroup. A character in the principal block has height zero
precisely when its degree is not divisible by $p$. Thus, in the
principal block the conjecture recovers one half of the
It\^o--Michler theorem: it characterizes when the Sylow $p$-subgroup
is abelian.

For many years it was believed that, when passing to blocks, one could
detect abelianity but not normality of the Sylow $p$-subgroup.
Surprisingly, in joint work with Schaeffer Fry \cite{ms},
we showed that blocks also detect normality, provided one considers
principal blocks for primes $q \neq p$.

\begin{theorem}
\label{thm:normal}
Let $G$ be a finite group, $p$ a prime, and $P$ a Sylow $p$-subgroup of
$G$. Then $P$ is normal in $G$ if and only if $p$ does not divide
$\chi(1)$ for every irreducible character $\chi$ lying in some
principal $q$-block, for some prime $q \neq p$.
\end{theorem}

Thus, by considering principal blocks for different primes, one recovers
exactly the normality condition that seemed to be lost from the
It\^o--Michler theorem. This result shows that blocks contain more local
information than was previously expected.

\subsection{The McKay conjecture}
\label{sub:mckay}

Among local-global problems, the one that has received the most attention
in recent decades is probably the McKay conjecture. Recall that for a
subgroup $H$ of $G$, the normalizer of $H$ is
\[
N_G(H) = \{ g \in G \mid H^g = H \}.
\]

Given a group $G$ and a prime $p$, I denote by $\Irr_{p'}(G)$ the
set of irreducible characters of $G$ whose degrees are not divisible by
$p$. The McKay conjecture, proposed in 1971 and recently proved by
Cabanes and Sp\"ath \cite{cs}, states the following.

\begin{conjecture}[McKay conjecture]
\label{conj:mckay}
Let $G$ be a finite group, $p$ a prime, and $P$ a Sylow $p$-subgroup of
$G$. Then
\[
|\Irr_{p'}(G)| = |\Irr_{p'}(N_G(P))|.
\]
\end{conjecture}

That is, the conjecture asserts that the number of irreducible characters
of degree not divisible by $p$ is the same in $G$ and in the local subgroup
$N_G(P)$.

In general, $N_G(P)$ is a much simpler group than $G$. For example, if
$p=2$ and $G$ is the symmetric group on $n$ letters, then the normalizer
of a Sylow $2$-subgroup is itself a Sylow $2$-subgroup. In this case, the
McKay conjecture asserts that the number of irreducible characters of odd
degree of the symmetric group equals $|P/P'|$.

Although it has been known for some time that this statement holds, even
in this case the underlying reason for the coincidence remains unclear.

Over the years, several stronger versions of the McKay conjecture have
been proposed, including the Alperin--McKay conjecture and refinements
due to Isaacs, Navarro, and Turull. Proofs of these stronger forms are
expected to be completed in the near future. However, it is not clear
that the current approaches are sufficient to explain why the McKay
conjecture is true, as the proofs rely on verifying complicated
conditions for all finite simple groups using the classification.

\subsection{The main problem of block theory}
\label{sub:main}

As mentioned in the introduction, in his article for the proceedings of
the Santa Cruz conference, Alperin \cite{alp} stated a number
of conjectures that have guided the subject for several decades. For
instance, \cite[Conj.~B]{alp} is the McKay conjecture,
\cite[Conj.~D]{alp} is Brauer's height zero conjecture,
and \cite[Conj.~E]{alp} is the Alperin--McKay conjecture.
Later on I will also state \cite[Conj.~C]{alp}, which
has now been proved and will play an important role in this article.

Here, however, I want to focus on another problem appearing at the very
beginning of Alperin's paper, one that does not seem to have received the
attention one might have expected. Alperin \cite{alp} begins
as follows:

\begin{quote}
Brauer’s work on block theory, stretching over decades, strongly suggests
the following problem as a reasonable choice for the main problem of the
subject.

Problem A. Give rules which determine the values of the characters of
$G$ in terms of the $p$-local subgroups of $G$.
\end{quote}

Alperin continues:

\begin{quote}
The values of the characters at the elements of order divisible by $p$ are
the most accessible but other values, for example degrees, have many
properties of a local nature. As a prime example of the sort of rule
desired we have McKay’s conjecture.
\end{quote}

This was in fact the second time that Alperin had singled out this
question as \emph{the main problem of block theory}; he had already done so
in \cite{alp2}.

As Alperin pointed out, the McKay conjecture is an excellent example of
the kind of rule envisioned in this program. However, it only concerns character values at the identity element,
or equivalently character degrees. The problem stated by Alperin is
much more ambitious: it asks for rules determining all character values in
terms of $p$-local subgroups.

Although Alperin explicitly presented this as ``the main problem of block
theory'', later developments largely focused on the usual local-global
conjectures, to the point that this aspect of his program seems to have
received far less attention than one might have expected. As far as I am
aware, the only place where this point is mentioned explicitly is Marcus
\cite{marcus}.

The statement of Alperin's problem is extraordinarily ambitious, and the
lack of progress for almost fifty years made it natural to suspect that it
might simply be out of reach.

The following sections suggest otherwise. As I shall explain, this problem
admits concrete and easily testable formulations for any given finite
group, and it can be approached from a surprisingly general perspective.

\section{Looking ahead}
\label{sec:future}

As indicated above, the conjectures and theorems in the previous section
essentially concern either abelian $p$-subgroups (Sylow subgroups or defect
groups) or characters of degree not divisible by $p$ (or, in the language
of blocks, of height zero). I have always regarded this height-zero and
abelian-defect perspective as only the first step. It seems natural to ask
whether the same local-global philosophy admits extensions to the case
where the Sylow $p$-subgroup or the defect group is nonabelian, or where
the prime $p$ is allowed to divide character degrees.

In this section I recall some conjectures, formulated in this direction
over the last two decades, that still concern character degrees and
heights. In the next section I turn to character values, where picky
elements and subnormalizers provide a more recent step in the same
program. In all cases there is significant evidence in their favor,
although it is reasonable to expect that a complete proof will require
many years.

\subsection{It\^o--Michler: nonabelian or nonnormal Sylow subgroups}
\label{sub:itonon}

I begin by discussing the possibility of obtaining versions of the
It\^o--Michler theorem in the case where the Sylow $p$-subgroup is not
normal or not abelian. I first introduce some notation.

For a finite group $G$, I define
\[
\cd(G)=\{\chi(1)\mid \chi\in\Irr(G)\},
\]
the set of irreducible character degrees of $G$.

A classical problem concerning this set, which will be relevant later on,
is the Isaacs--Seitz conjecture.

\begin{conjecture}[Isaacs--Seitz]
\label{conj:is}
If $G$ is solvable, then
\[
\dl(G)\le |\cd(G)|.
\]
\end{conjecture}

This conjecture was formulated in the 1960s. It was solved in many cases
during the 1980s and 1990s, but the general problem remains open. The best
general result from that period, due to Gluck, shows that
\[
\dl(G)\le 2|\cd(G)|.
\]

Huppert's work on certain families of $p$-groups suggested that one should
expect a logarithmic bound, which led to substantial activity in the area.
This eventually culminated in a series of papers by Keller, where he
obtained the bound
\[
\dl(G/F(G))\le 24\log_2 |\cd(G)|+364.
\]
Somewhat paradoxically, this reduced the search for logarithmic bounds in
the solvable case to the case of $p$-groups, precisely the case from which
the original conjecture had emerged.

I now return to local-global problems. Fix a prime $p$, and define
\[
\cd_p(G)=\{\chi(1)_p\mid \chi\in\Irr(G)\},
\]
the set of $p$-parts of the irreducible character degrees of $G$.

In this notation, the It\^o--Michler theorem says that, if $P$ is a Sylow
$p$-subgroup of $G$, then
\[
\dl(P)=1 \quad\text{and}\quad N_G(P)=G
\]
if and only if
\[
\cd_p(G)=\{1\}.
\]

This suggests that $\cd_p(G)$ might control the derived length of
$P$. Since $p$-groups are solvable, this is a natural idea. This led to
the following conjectures, formulated in \cite{mor03} and
announced at the conference celebrating Thompson's 70th birthday in
Gainesville in 2003.

\begin{conjecture}
\label{conj:cdp}
Let $G$ be a finite group, $p$ a prime, and $P$ a Sylow $p$-subgroup of
$G$. Let
\[
p^b=\max \cd(P)
\qquad\text{and}\qquad
p^f=\max \cd_p(G).
\]
Then:
\begin{enumerate}
\item $\dl(P)$ is bounded in terms of $f$;
\item $b$ is bounded in terms of $f$;
\item $\dl(P)$ is bounded in terms of $|\cd_p(G)|$.
\end{enumerate}
Moreover, the bounds in (1) and (3) should be logarithmic.
\end{conjecture}

It was further suggested that one might even have $b\le 2f$. This bound
would imply the logarithmic bound in (1), and (1) also follows from (3).

A few years later, in joint work with Isaacs, Navarro, and Tiep, I
formulated the following stronger conjecture.

\begin{conjecture}
\label{conj:imnt}
Let $G$ be a finite group, $p$ a prime, and $P$ a Sylow $p$-subgroup of
$G$. Then
\[
|\cd(P)|\le |\cd_p(G)|+1.
\]
\end{conjecture}

The It\^o--Michler theorem implies that this conjecture holds in the case
where $|\cd_p(G)|=1$. In \cite{imnt} it was proved in the case where
$|\cd_p(G)|=2$.

This conjecture implies Conjecture~\ref{conj:cdp}(3), and, assuming a
logarithmic bound in the Isaacs--Seitz problem, it would also give the
expected logarithmic form of that bound.

These conjectures may be viewed as natural extensions of the abelian part
of the It\^o--Michler theorem to the case of nonabelian Sylow
$p$-subgroups.

There is another fundamental aspect of It\^o--Michler, namely the
normality of the Sylow $p$-subgroup. In \cite{gglmnt}
we observed that this part might also admit an extension if one considers
the number of irreducible characters of degree divisible by $p$. In
particular, I proved that if $G$ has exactly one irreducible character of
degree divisible by $p$, then $N_G(P)$ is a maximal subgroup of $G$.
Recall that $P$ is normal if and only if $N_G(P)=G$. Thus, having ``almost
no'' irreducible characters of degree divisible by $p$ forces $N_G(P)$ to
be ``almost'' equal to $G$.

This led to the following conjecture.

\begin{conjecture}
\label{conj:max}
Let $G$ be a finite group, $p$ a prime, and $P$ a Sylow $p$-subgroup of
$G$. If $G$ has $n$ irreducible characters of degree divisible by $p$,
then the length of any subgroup chain between $N_G(P)$ and $G$ is at most
$n$.
\end{conjecture}

For example, if $N_G(P)$ is maximal, the only subgroup chain between
$N_G(P)$ and $G$ is
\[
N_G(P)<G,
\]
which has length $1$.

\subsection{Brauer's height zero conjecture: nonabelian defect}
\label{sub:bhznon}

Let $B$ be a $p$-block of $G$. I write $\operatorname{ht}(B)$ for the set of
heights of the irreducible characters in $B$. Brauer's height zero
conjecture asserts that $\dl(D)=1$, where $D$ is the defect group,
if and only if $|\operatorname{ht}(B)|=1$ (Brauer proved that every block has
characters of height zero).

Can one say something when $|\operatorname{ht}(B)|>1$? In \cite{mor03}
I proposed the following partial extensions of Brauer's height zero
conjecture.

\begin{conjecture}
\label{conj:ht}
Let $G$ be a finite group, $p$ a prime, and $B$ a $p$-block of $G$ with
defect group $D$. Let
\[
p^b=\max \cd(D)
\qquad\text{and}\qquad
e=\max \operatorname{ht}(B).
\]
Then:
\begin{enumerate}
\item $\dl(D)$ is bounded in terms of $e$;
\item $b$ is bounded in terms of $e$;
\item $\dl(D)$ is bounded in terms of $|\operatorname{ht}(B)|$.
\end{enumerate}
Moreover, the bounds in (1) and (3) should be logarithmic.
\end{conjecture}

In \cite{mor03} I proved these conjectures, as well as
Conjecture~\ref{conj:cdp}, for general linear groups and symmetric groups (in the latter
case, only the details for $p>3$ appeared there, although the complete
proof was later published in \cite{gms}).

I should point out that these conjectures, as well as Conjecture~\ref{conj:cdp}, have
not been widely explored so far. In fact, related conjectures of more
limited scope have appeared in some recent papers, without explicit
reference to these earlier conjectures.

These inequalities seem to capture structural complexity in a fundamental
way within the general local-global framework, and may be of the same
interest as the problems discussed in the previous section. The most prominent example is Brauer's $k(B)$-conjecture,
the only one among the classical local-global problems that has not been
reduced to finite simple groups, despite having been studied extensively.
This conjecture asserts that the number of irreducible characters in a
block $B$ does not exceed the order of its defect group.

I also proposed the following additional inequality, which I consider
natural.

\begin{conjecture}
\label{conj:cdht}
Let $G$ be a finite group, $p$ a prime, and $B$ a $p$-block of $G$ with
defect group $D$. Then
\[
|\cd(D)|\le |\operatorname{ht}(B)|+1.
\]
\end{conjecture}

This is, of course, the block version of Conjecture~\ref{conj:imnt}.

It is interesting to compare the set of character degrees of $D$ with the
set of heights of the characters in $B$. In general, one cannot expect any
strong relation between these two sets; inequalities such as the one above
seem to be the most one can reasonably hope for at this level of
generality.

However, there is one remarkable feature that these two sets do appear to
share: the smallest nontrivial element.

\begin{conjecture}
\label{conj:em}
Let $G$ be a finite group, $p$ a prime, and $B$ a $p$-block of $G$ with
defect group $D$. Then
\[
\inf(\cd(D)\setminus\{1\})=p^{\inf(\operatorname{ht}(B)\setminus\{0\})}.
\]
\end{conjecture}

By convention, I take the infimum of the empty set to be infinity. Note
that $\cd(D)\setminus\{1\}$ is empty if and only if $D$ is
abelian, and $\operatorname{ht}(B)\setminus\{0\}$ is empty if and only if all
characters in $B$ have height zero. Thus, the fact that both sides are
simultaneously infinite is exactly Brauer's height zero conjecture.

This generalization was proposed jointly with Charles Eaton in 2014
\cite{em}, although it had already been announced at the
conference in honor of Isaacs held in Valencia in 2009.

After some initial skepticism, confidence in the conjecture has gradually
increased. Recent progress has been significant, although I am still far
from a complete proof. In \cite{em} it was shown that the
inequality
\[
\inf(\cd(D)\setminus\{1\})\le p^{\inf(\operatorname{ht}(B)\setminus\{0\})}
\]
follows from Dade's projective conjecture, another central problem in the
theory. This conjecture has seen major progress, and its reduction to
finite simple groups has been established.

The more delicate part is the reverse inequality,
\[
\inf(\cd(D)\setminus\{1\})\ge p^{\inf(\operatorname{ht}(B)\setminus\{0\})}.
\]
This has recently been proved in several important cases, including:
principal blocks whose defect group has two character degrees, in
\cite{mmr}; for $p\ge 5$, there are no minimal counterexamples
among finite simple groups \cite{ms1}; and principal
blocks of $p$-solvable groups \cite{nav24}.

Now that Brauer's height zero conjecture has been proved, it is reasonable
to expect that interest in this generalization will increase considerably.

\section{Picky elements and subnormalizers}
\label{sec:picky}

\subsection{Picky elements and a first step toward Alperin's Problem~A}
\label{sub:picky}

Sometimes the development of a mathematical idea begins with a phenomenon
that at first sight appears rather minor. This was exactly the case for
the results described in this subsection.

In 2022 I began working with Attila Mar\'oti (Alfr\'ed R\'enyi Institute,
Budapest) and with my then PhD student Juan Mart\'inez Madrid on an
apparently elementary problem concerning the covering of the set of
$p$-elements of a finite group by its Sylow $p$-subgroups. Recall that,
by Sylow's theorem, every element of $p$-power order belongs to some Sylow
$p$-subgroup. In other words, the Sylow $p$-subgroups of $G$ cover the set
of $p$-elements of $G$.

The question that I proposed, and which initiated our work, was easy to state:

\begin{quote}
For which finite groups $G$ are all Sylow $p$-subgroups needed in order to
cover the set of $p$-elements?
\end{quote}

The key step was the following elementary lemma.

\begin{lemma}
\label{lem:cover}
Let $G$ be a finite group and let $p$ be a prime. Then all Sylow
$p$-subgroups are needed in order to cover the set of $p$-elements of $G$
if and only if there exists a $p$-element $x\in G$ that belongs to a
unique Sylow $p$-subgroup.
\end{lemma}

Somewhat unexpectedly, in \cite{mmm} we showed that most
finite groups do contain such elements. In other words, it is typical
rather than exceptional that every Sylow $p$-subgroup be needed in order
to cover the set of $p$-elements.

It soon became clear that these elements have strong structural properties
from the point of view of character theory. In fact, we proved the
following.

\begin{proposition}
\label{prop:zero}
Let $x\in G$ be a $p$-element lying in a unique Sylow $p$-subgroup. Then
\[
\chi(x)=0
\]
for every $\chi\in\Irr(G)$ that does not belong to a block whose defect
group is a Sylow $p$-subgroup of $G$.
\end{proposition}

The behavior of zeros of irreducible characters has long been an important
source of structural information. Recall, for instance, Burnside's
classical theorem that an irreducible character is nonlinear if and only
if it has a zero; see Isaacs~\cite{isa}.

On June 22nd, 2023 I gave a talk on \cite{mmm} at a conference
in honor of Pavel Shumyatsky. The following days I began to wonder whether
the relevance of these elements in character theory might be related to
Alperin's Conjecture~C.

For a subset $S\subseteq G$, write
\[
\Irr^S(G)=\{\chi\in\Irr(G)\mid \chi(s)\neq 0
\text{ for some }s\in S\}.
\]

\begin{conjecture}[Alperin's Conjecture~C]
\label{conj:alp}
Let $G$ be a finite group, let $p$ be a prime, and let $P$ be a Sylow
$p$-subgroup of $G$. Assume that $P$ is TI, that is,
\[
P^g\cap P\in\{1,P\}\qquad\text{for all }g\in G.
\]
Then
\[
|\Irr^P(G)|=|\Irr(N_G(P))|.
\]
\end{conjecture}

As Alperin observed, one has
\[
\Irr(N_G(P))=\Irr^P(N_G(P)),
\]
so the conjecture says that when $P$ is TI, the number of irreducible
characters that do not vanish identically on $P$ is the same in $G$ and
in the local subgroup $N_G(P)$. This conjecture was proved by Blau and
Michler~\cite{bm} in 1990 using the classification of the finite
simple groups.

The connection with the elements considered above is immediate:
\[
P \text{ is TI } \iff \text{ every nontrivial element of } P
\text{ is picky.}
\]

At that point I realized that this phenomenon might admit a much more general
formulation. To explore this possibility, I turned to GAP~\cite{gap}. The
computational evidence was strikingly robust, and strongly suggested that
one was not looking at isolated coincidences. This led me to formulate the
following conjecture.

Let $x\in P$, and write
\[
\Irr^x(G)=\{\chi\in\Irr(G)\mid \chi(x)\neq 0\}.
\]

\begin{conjecture}
\label{conj:picky}
Let $G$ be a finite group, let $p$ be a prime, and let $x$ be a $p$-element
lying in a unique Sylow $p$-subgroup. Then there exists a bijection
\[
f:\Irr^x(G)\longrightarrow \Irr^x(N_G(P))
\]
such that:
\begin{enumerate}
\item $\chi(1)_p=f(\chi)(1)_p$ for all $\chi\in\Irr^x(G)$;
\item $\mathbb{Q}(\chi(x))=\mathbb{Q}(f(\chi)(x))$ for all
      $\chi\in\Irr^x(G)$.
\end{enumerate}
\end{conjecture}

In July 2023 I communicated these ideas to Noelia Rizo. She suggested the
name \emph{picky elements} for the $p$-elements lying in a unique Sylow
$p$-subgroup. We then began working on Conjecture~\ref{conj:picky}, which
we will refer to as the \emph{Picky Conjecture}, and later on the
subnormalizer conjecture \cite{mr}.

In many cases one even has the stronger property
\[
\chi(x)=\pm f(\chi)(x).
\]
In that case, I say that $(G,x)$ satisfies the \emph{strong picky
conjecture}.

Among the cases that have already been proved are the following:
\begin{itemize}
\item the strong picky conjecture holds for symmetric groups
      \cite{mar};
\item the strong picky conjecture holds for finite groups of Lie type in
      characteristic $q\neq p$ \cite{ms2};
\item the strong picky conjecture holds for $p$-solvable groups when
      $p>2$ \cite{mnr}.
\end{itemize}

The main examples where the strong picky conjecture fails are precisely
certain picky elements of finite simple groups with nonabelian TI Sylow
subgroups. In those examples one has instead the following relation
between the values of $\chi$ and $f(\chi)$ at $x$, suggested by
Gunter Malle:
\[
\chi(x)_p=f(\chi)(x)_p
\qquad\text{for all }\chi\in\Irr^x(G).
\]

Here, as in Isaacs~\cite{isa}, the values of irreducible characters
are algebraic integers, and the $p$-part of an algebraic integer $\alpha$
may be defined by
\[
\alpha_p=
\left|N_{\mathbb{Q}(\alpha)/\mathbb{Q}}(\alpha)\right|_p^{1/[\mathbb{Q}(\alpha):\mathbb{Q}]}.
\]

Thus, the picky conjecture relates the values of the irreducible
characters of a group $G$ at a picky element $x$ with the values at the
same element of the irreducible characters of the normalizer of a Sylow
$p$-subgroup. A proof of this conjecture would therefore provide a local rule for the
character values at picky elements, in the spirit of Alperin's ``main
problem of block theory''.

Moreover, condition~(1) in the picky conjecture tells us that the bijection
preserves the $p$-parts of character degrees. This is also significant. A
basic result in character theory asserts that an irreducible character of
degree not divisible by $p$ never vanishes on a $p$-element. Thus
\[
\Irr_{p'}(G)\subseteq \Irr^x(G)
\]
for every $p$-element $x$. Since the bijection in the picky conjecture
preserves $p$-parts of degrees, it sends the $p'$-degree characters of
$G$ to the $p'$-degree characters of $N_G(P)$. In other words, the picky
conjecture implies the McKay conjecture for groups with picky elements.

The versions stated here are the most elementary ones, for simplicity.
More general forms will appear in the forthcoming paper~\cite{mor}.

\subsection{The subnormalizer conjecture}
\label{sub:sub}

What happens in connection with Alperin's Problem~A for elements that are
not picky? The surprises did not end with the picky conjecture; rather,
that was only the beginning.

My work in \cite{mmm} led me to consider the
subnormalizer, studied in a series of papers by Carlo Casolo
\cite{cas1,cas2,cas3}. Whereas the normalizer of a
subgroup is a completely standard notion in group theory, the
subnormalizer has so far received less attention than it deserves.

Recall that a subgroup $H$ of a group $G$ is called subnormal in $G$ if
there exists a chain of subgroups
\[
H=H_0\leq H_1\leq \cdots \leq H_t=G
\]
such that $H_i\trianglelefteq H_{i+1}$ for every $i=0,\dots,t-1$. In this
case I write $H\trianglelefteq\trianglelefteq G$.

Given any subgroup $H$ of $G$, the subnormalizer of $H$ in $G$ is defined
by
\[
S_G(H)=\{\,g\in G \mid H\trianglelefteq\trianglelefteq \langle H,g\rangle\,\}.
\]
In general, $S_G(H)$ is not a subgroup of $G$, but merely a subset.
Accordingly, I define the \emph{subnormalizer subgroup} by
\[
\Sub_G(H)=\langle S_G(H)\rangle.
\]

In this article I will only need to consider subnormalizer subgroups of
cyclic subgroups. If $H=\langle x\rangle$, I write
\[
\Sub_G(x)=\Sub_G(H).
\]

The key observation leading to the subnormalizer conjecture is the
following consequence of Casolo's work.

\begin{lemma}
\label{lem:sub}
Let $G$ be a finite group, let $p$ be a prime, let $P$ be a Sylow
$p$-subgroup of $G$, and let $x\in P$. Then:
\begin{enumerate}
\item $N_G(P)\subseteq \Sub_G(x)$;
\item $N_G(P)=\Sub_G(x)$ if and only if $x$ lies in a unique Sylow
      $p$-subgroup.
\end{enumerate}
\end{lemma}

This observation marked a turning point. As soon as I saw part~(2), I
immediately thought that I had been overlooking something important. If,
as the picky conjecture suggests, the set $\Irr^x(G)$ captures more
faithfully than $\Irr_{p'}(G)$ the information relevant to a given
$p$-element $x$, then perhaps the subgroup that really encodes the
character values of $G$ at $x$ is not $N_G(P)$, but rather $\Sub_G(x)$.

Computations in GAP immediately confirmed that something deep was going
on. For instance, if $x$ is an $8$-cycle in the symmetric group $S_{16}$,
then the list of nonzero character values (up to sign), together with the
$2$-part of the degree and their multiplicities, is the same in $S_{16}$
and in
\[
\Sub_{S_{16}}(x)=S_8\wr C_2,
\]
as shown in Table~\ref{tab:s16}.

\begin{table}[H]
\centering
\small
\renewcommand{\arraystretch}{0.95}
\begin{tabular}{|c|c|c||c|c|c|}
\hline
Value & $2$-part & Multiplicity & Value & $2$-part & Multiplicity \\
\hline
$1$  & $1$   & $4$ & $36$ & $2$   & $1$ \\
$7$  & $1$   & $2$ & $14$ & $4$   & $8$ \\
$35$ & $1$   & $2$ & $42$ & $4$   & $5$ \\
$6$  & $2$   & $1$ & $70$ & $4$   & $6$ \\
$14$ & $2$   & $2$ & $90$ & $4$   & $2$ \\
$20$ & $2$   & $1$ & $20$ & $8$   & $6$ \\
$28$ & $2$   & $2$ & $28$ & $8$   & $7$ \\
$34$ & $2$   & $1$ & $56$ & $16$  & $6$ \\
     &       &     & $64$ & $128$ & $9$ \\
\hline
\end{tabular}
\caption{Absolute character values of $S_{16}$ on an $8$-cycle, the
$2$-part of the degree, and their multiplicities.}
\label{tab:s16}
\end{table}

This led to the following conjecture.

\begin{conjecture}[The Subnormalizer Conjecture]
\label{conj:sub}
Let $G$ be a finite group, let $p$ be a prime, and let $x$ be a
$p$-element. Then there exists a bijection
\[
f:\Irr^x(G)\longrightarrow \Irr^x(\Sub_G(x))
\]
such that:
\begin{enumerate}
\item $\chi(1)_p=f(\chi)(1)_p$ for all $\chi\in\Irr^x(G)$;
\item $\mathbb{Q}(\chi(x))=\mathbb{Q}(f(\chi)(x))$ for all
      $\chi\in\Irr^x(G)$.
\end{enumerate}
\end{conjecture}

By Lemma~\ref{lem:sub}, the picky conjecture is precisely the special case
of this conjecture in which $x$ lies in a unique Sylow $p$-subgroup. As
in the picky conjecture, I expect that in many cases the stronger
relation
\[
\chi(x)=\pm f(\chi)(x)
\]
should hold for all $\chi\in\Irr^x(G)$. This has already been proved, for
instance, for symmetric groups when $p=2$ in work of Mart\'inez Madrid
\cite{mar}, and for several finite groups of Lie type in
work of Malle \cite{mal1,mal2}.

Thus, the subnormalizer conjecture proposes what had remained elusive for
fifty years: for a $p$-element $x$, the local subgroup in which one should
look for information on the values of the irreducible characters of $G$ at
$x$ is $\Sub_G(x)$.

Here I have only formulated the conjecture for $p$-elements, but a variant
of it makes sense for essentially arbitrary elements of $G$. Whereas the
subnormalizers of $p$-elements tend to be large subgroups (for example,
they always contain $N_G(P)$), the subnormalizers of elements whose order
is divisible by more than one prime tend to be surprisingly small. Thus,
the conjecture suggests that the character values of $G$ at mixed-order
elements may be encoded in subgroups of unexpectedly small size. Even
after working extensively with examples, I still find this remarkable.

An additional difficulty is that, although several properties of
subnormalizers of elements of prime-power order are known, the structure of the
subnormalizers of mixed-order elements remains largely mysterious. In
particular, a substantial part of the work of Mart\'inez Madrid for
symmetric groups and of Malle for groups of Lie type consists precisely in
understanding the structure of $\Sub_G(x)$ in these cases.

\section{Conclusion}
\label{sec:concl}

Just as the classification of the finite simple groups was only the
starting point for a systematic attack on local-global problems in
representation theory, the recent advances in this area should be viewed
as opening a new stage rather than closing the subject.

The McKay conjecture has now been proved by means of a very long argument,
one that has already produced many new ideas with impact on other
problems. However, it is still not understood why the conjecture is true.

The subnormalizer conjecture proposes two significant extensions of the
classical McKay framework. On the one hand, it replaces the set
$\Irr_{p'}(G)$ by the richer set $\Irr^x(G)$. On the other hand, it
replaces the subgroup $N_G(P)$ by the subgroup $\Sub_G(x)$. With these
two structural changes, everything seems to work in a more natural and
coherent way. For instance, as mentioned in Subsection~\ref{sub:bhznon}, within the
classical framework the connections between local and global information
are quite limited: even comparatively modest conjectures such as
Conjecture~\ref{conj:em} were initially met with skepticism. That conjecture only
predicts a relation between the smallest positive height occurring in a
block and the smallest nontrivial character degree of its defect group.
By contrast, the subnormalizer conjecture suggests something much
stronger: that once one chooses the ``correct'' sets of characters,
namely $\Irr^x(G)$ and $\Irr^x(\Sub_G(x))$, there should exist a
bijection preserving all $p$-parts of character degrees.

As indicated in \cite{mor25}, this perspective goes well
beyond the picky and subnormalizer conjectures stated here. Quite
naturally, one is led to try to reformulate in terms of $\Irr^x(G)$ any
problem that has traditionally been posed for $\Irr_{p'}(G)$. Several
examples are mentioned in \cite{mor25}: conjectures that
are false in the classical setting, but are very likely to hold from this
new point of view.

The same source also points out connections with Brou\'e's famous
conjecture for abelian defect groups, one of the major open problems in
representation theory. It is known that Brou\'e's conjecture fails for
nonabelian defect groups, but the more general results to appear in the forthcoming paper~\cite{mor}
suggest that meaningful extensions may exist beyond the abelian case.

There are many other connections that lie beyond the scope of the present
article. At this point, I do not know which of them, if any, will
eventually provide the conceptual explanation for what is happening in
these local-global conjectures (or theorems). It does seem, however, that
some missing concept is still needed, perhaps not purely group-theoretic,
and not necessarily algebraic, in order to unify these ideas.

For instance, it is striking that simplicial complexes have appeared in
the study of several local-global conjectures. Years before these later
connections, Casolo \cite{cas3} had already used simplicial
complexes in his study of the subnormalizer, in a context in which
representations played no role at all. That work did not influence the
later appearance of simplicial complexes in local-global problems, and
yet the coincidence is hard to ignore.

The work in \cite{mmm} and the more general results to appear in
\cite{mor} also suggest that fusion systems may be related to these
conjectures. In
\cite{mmm} we wrote: ``It seems reasonable to think that
fusion systems could be helpful to study the condition provided by
Lemma~2.1.'' The lemma referred to there is exactly Lemma~\ref{lem:cover}, so
already in \cite{mmm} we were suggesting that fusion
systems might be useful in the study of picky elements.

One of the results in the forthcoming paper~\cite{mor} that convinced me
that the subnormalizer is the right subgroup in which to look for local
information on the character values at a given element is the following.

\begin{lemma}
\label{lem:fusion}
Let $G$ be a finite group, let $p$ be a prime, let $P$ be a Sylow
$p$-subgroup of $G$, and let $x\in P$. If $y\in \Sub_G(x)$ is conjugate
to $x$ in $G$, then $x$ and $y$ are conjugate in $\Sub_G(x)$.
\end{lemma}

Although this lemma can also be obtained using some of Casolo's results
on the subnormalizer, my original proof used Alperin's fusion theorem,
which in turn was one of the origins of the theory of fusion systems. For
an introduction to fusion systems, I recommend Broto's article
\cite{bro} in this same section of \emph{La Gaceta}.

These conjectures have already generated interest in different contexts.
For example, they appear in connection with the work of Giudici, Morgan,
and Praeger \cite{gmp} on permutation groups. All this
suggests that the ideas presented in this article will lead to active
lines of research for a long time.


\end{document}